\documentclass[12pt]{article}

\usepackage{amsmath}
\usepackage{amsthm}
\usepackage{color}

\def\uo{\ensuremath{\underline{0}}}
\usepackage[colorlinks=true,
linkcolor=green,
filecolor=brown,
citecolor=green]{hyperref}

\begin{document}

\begin{center}
    \textbf{A Combinatorial Proof of Sun's ``Curious'' Identity}

\vskip 20pt

{\bf David Callan}\\
 Department of Statistics,  
University of Wisconsin-Madison \\ 
1210 W. Dayton Street,  
Madison, WI 53706-1693\\
\texttt{callan@stat.wisc.edu} 

\vskip 5pt

January 17 2004

\end{center}

\vskip 5pt

\centerline{\bf Abstract}
\noindent
{A binomial coefficient identity due to Zhi-Wei 
Sun is the subject of half a dozen recent papers that prove it by 
various analytic techniques and establish a generalization. Here we 
give a simple proof that uses weight-reversing involutions on suitable 
configurations involving dominos and colorings.
With somewhat more work, the method extends to the generalization also.
 }\\

\vskip 30pt

{\Large \textbf{0 \quad Introduction}  }

The identity
\begin{equation}
\sum_{i=0}^{m}(x+m+1) {(-1)^i {{x+y+i} \choose {m-i}} {{y+2i} \choose
{i}}}
-\sum_{i=0}^{m} {{{x+i} \choose {m-i}}(-4)^i}=(x-m) {{x} \choose {m}}
\label{eq:1}
\end{equation}
was proved \cite{sun,panprod,merlini,ekhad,jensen} using induction, generating functions, Riordan arrays, the 
WZ method, and the so-called Jensen formula respectively.
This identity is the case $z=1$ of the following generalization 
\cite{extension}:
\begin{align} \label{eq:2}
   &   
   \left(x+(m+1)z\right)\sum_{n=0}^{m}(-1)^n\binom{x+y+nz}{m-n}\binom{y+n(z+1)}{n}  \\ \notag
 = & \ z\sum_{0\le l\le n \le m}(-1)^n\binom{n}{l}\binom{x+l}{m-n}(1+z)^{n+l}(1-z)^{n-l}+(x-m)\binom{x}{m}.
\end{align}
Section 1 is a ``lite'' version of the paper: a proof of (\ref{eq:1}) using 
weight-reversing involutions on configurations on an initial segment 
of the positive integers that involve
some vertices covered by dominos and the rest colored 
black or white. Section 2 then uses the same method in a somewhat more 
complicated setting to prove (\ref{eq:2}). 

{\Large \textbf{1 \quad Identity (1)}  }

Both sides of (\ref{eq:1}) are polynomials in $x$ of degree $m+1$ 
with the same leading coefficient, $1/m!$, and so 
it suffices to prove (\ref{eq:1}) for $m+1$ distinct values 
of $x$. The right hand side vanishes for $x=0,1,2,\ldots,m$ and so, 
putting $x=m-k$, replacing the summation index $i$ by $m-i$ and 
canceling a $(-1)^{m}$ factor, (\ref{eq:1}) is equivalent to
\begin{multline*}
    (2m-k+1)\sum_{i=0}^{m}(-1)^{i}\binom{2m+y-k-i}{i}\binom{2m+y-2i}{m-i} = \\
    \sum_{i=0}^{m}(-1)^{i}\binom{2m-k-i}{i}2^{2(m-i)} \hspace*{10mm} 
    0 \le k \le m.
\end{multline*}
We show that in fact both sides $=(2m-k+1)2^{k}$. 
 
After some cancellation, this amounts to
\begin{align}
\sum_{i=0}^{m}(-1)^{i}\binom{2m+y-k-i}{i}\binom{2m+y-2i}{m-i} \quad & = \quad
2^{k},\hspace*{10mm} \label{eq:3}\\ \nonumber
 \textrm{and}\hspace*{10mm} & & \\
 \sum_{i=0}^{m}(-1)^{i}\binom{2m-k-i}{i}2^{2m-2i} \quad & = \quad 
 (2m-k+1)2^{k} \label{eq:4}
\end{align}

Fix an integer $k \in [0,m]$, set $y=b$, a fixed nonnegative 
integer and let $[n]$ denote $\{1,2,\ldots,n\}$.  
The right hand side of (\ref{eq:3}) is the total weight of what we will call 
Sun(3)-configurations on $[2m+b]$ defined as follows: cover 
some of $[2m+b-k]$ with $i$ ($0\le i \le m$) nonoverlapping dominos, each covering two 
adjacent vertices [$\binom{2m+b-k-i}{i}$ choices] and then color $m-i$ of the remaining 
$2m+b-2i$ vertices black and the rest white [$\binom{2m+b-2i}{m-i}$ choices].
Assign weight $=(-1)^{\# \textrm{dominos}}$.
Two examples are given in the Figure below (both for $m=4,\,b=2,\,k=1$ and so 
$2m+b-k=9$; an underline 
denotes a black vertex). The first has weight $-1$, the second has 
weight $+1$.
\[
\underbrace{1\ \underline{2}\  \underline{3}\ 4\ 5\ 6\ \fbox{7}\fbox{8}\  
9}_{\textrm{domino range}}\ \underline{10} \hspace*{15mm} 
\underbrace{1\ \underline{2}\  \fbox{3}\fbox{4}\ 5\ 6\ \fbox{7}\fbox{8}\  
\ 9}_{\textrm{domino range}}\ \underline{10}
\]
The left hand side of (\ref{eq:3}) is the total weight of all Sun(3)-configurations. 
Here is a weight-reversing involution on most of them.
Call $[2m+b-k]$ the \emph{domino 
range} and say two consecutive vertices form an \emph{active pair} if both lie 
in the domino range and are either (i) covered by a domino or (ii) colored 
black, white in that order. 
Look for the leftmost active pair of vertices. If active by virtue of 
(i), remove the domino and color the left vertex black, the right one 
white. If active by virtue of 
(ii), remove the colors and cover them with a domino. For example, 
this map 
interchanges the two configurations above. The map 
preserves \#~dominos + \# black vertices while altering \# dominos by 
1 and hence is a 
weight-reversing involution on all Sun(3)-configurations except
those with no active pairs. Such an exceptional configuration entails no dominos ($i=0$) and hence 
weight $=1$, and all black vertices in the domino range flush left. 
When $i=0$, the vertices outside the domino range may be colored arbitrarily 
($j\in[0,k]$ black vertices outside implies $m-k+b+j$ white vertices 
inside and $m-k+b+j\ge 0$ since $k\le m$) and the flush left 
requirement then determines the configuration. Hence $2^{k}$ exceptional 
configurations, and (\ref{eq:3}) follows.

Identity (\ref{eq:4}) can be proved similarly using the following 
configurations: place $i\ge 0$ dominos 
on $[2m-k]$ [$\binom{2m-k-i}{i}$ choices] but now color the remaining 
$2m-2i$ vertices black or white independently [$2^{2m-2i}$ choices].
Assign weight $=(-1)^{\# \textrm{dominos}}$ and use the same map. 
Again, the exceptional configurations are those with no dominos and all black 
vertices in the domino range flush left.
The number of such 
configurations is now $(2m-k+1)$\,[locate switchover from black to white 
in domino range] $\times 2^{k}$\,[color vertices outside domino 
range], and (\ref{eq:4}) follows.

\vspace*{10mm}

{\Large \textbf{2 \quad Identity (2)}  }

Just as for (\ref{eq:1}) and with the same change of variable, (\ref{eq:2})
is equivalent to
\begin{multline*}   
   \left( (m-k)+(m+1)z\right)\sum_{i=0}^{m}(-1)^{(m-i)}\binom{m-k+y+(m-i)z}{i}
   \binom{y+(m-i)(z+1)}{m-i}  \\ 
 = \ z\sum_{0\le j\le i \le m}(-1)^i\binom{i}{j}
 \binom{m-k+j}{m-i}(1+z)^{i+j}(1-z)^{i-j}
 \hspace*{20mm} 0 \le k \le m
\end{multline*}
Again, both sums have a closed form and it suffices to show that with 
$y=b$ a nonnegative integer and $z=q$ also a nonnegative integer,
\begin{multline}
     \sum_{i=0}^{m}(-1)^{i}\binom{(q+1)m-k+b-qi}{i}
     \binom{(q+1)m+b-(q+1)i}{m-i}
     =q^{m-k}(1+q)^{k}
    \label{omino}
\end{multline}
and
\begin{multline}
    \sum_{0 \le j \le i \le  m}^{}(-1)^{m-i}\binom{i}{j}
    \binom{m-k+j}{m-i}(1+q)^{i+j-k}(1-q)^{i-j}\\
    =\ \big( (m-k)+(m+1)q\big)q^{m-k-1}\hspace*{20mm}
    \label{matrix}
\end{multline}
both for $ 0 \le k \le m$.

For (\ref{omino}), the analog of a Sun(3)-configuration uses 
$(q+1)$-ominos i.e., $q+1$ consecutive vertices: choose $i \le m\ (q+1)$-ominos 
in the ``omino range'' $[(q+1)m+b-k]$, choose $m-i$ of the remaining 
$(q+1)m+b-(q+1)i$ vertices in $[(q+1)m+b]$ to color black ($B$) and the rest 
white ($W$). Use weight $=(-1)^{\mathrm{\# ominos}}$. Here the 
active $(q+1)$-tuples in the omino range  are  $(q+1)$-ominos and 
strings $W^{q}B\ (=\underbrace{W\ldots W}_{q}B)$. The analogous 
involution is obvious and the survivors are again the configurations 
with no active $(q+1)$-tuple. This entails no omino (so $i=0$ and 
weight $=1$) and no $W^{q}B$ in $[(q+1)m+b-k]$. Such a configuration 
has arbitrary coloring of the last $k$ vertices, say $j$ $B$s with 
$0\le j \le k$ and then begins $W^{i_{1}}B W^{i_{2}}B \ldots 
W^{i_{m-j}}BW\ldots W$ with each exponent $i_{\ell} \in [0,q-1]$. We find 
that there are $\sum_{j=0}^{k}\binom{k}{j}$\,[choose $B$s among last $k$ 
vertices]\,$\times q^{m-j}$\,[choose exponents]$ =q^{m-k}(1+q)^{k}$ 
survivors, as expected.

The form of the right hand side of (\ref{matrix}) suggests separate 
treatment of the case $k=m$:
 \[
\sum_{0\le j \le 
i\le m}(-1)^{m-i}\binom{i}{j}\binom{j}{m-i}(1+q)^{i+j-m}(1-q)^{i-j}=m+1
\]
Here, the relevant configurations are 5-tuples $(i,J,K,A,B)$ with $i \in 
[0,m],\ J$ a $j$-subset of $I:=[i],\ K$ an $(m-i)$-subset of $J,\ A$ 
an arbitrary $q$-colored subset of $J\backslash K$ and $B$ 
an arbitrary $q$-colored subset of $I\backslash J$. The weight is 
$(-1)^{m-i+\vert B \vert}$. Since the value of $i$ is built into $K$, 
whose size is $m-i$, and the conditions $J,B \subseteq I$ say that 
$\max\{J,B\}\le m -\vert K \vert$, these configurations may be more 
compactly described as 4-tuples $(J,K,A,B)$ satisfying $J \subseteq 
[m]$, both $K$ and $A \subseteq J,\ K \cap A =\emptyset,\ J \cap B = 
\emptyset,\ A$ and $B$ both $q$-colored, and finally, $\max\{J,B\}\le m -\vert K \vert$.
The weight is $(-1)^{\vert  K \vert+\vert B\vert}$. Such a 4-tuple can be 
conveniently represented by  a $2 \times m$ 0-1 matrix with (possibly) 
some underlined 0s: 1s in the top (resp. bottom) row indicate 
elements of $K$ (resp. $J$) and underlined 0s in the top (resp. bottom) row indicate 
elements of $A$ (resp. $B$). The weight is then $(-1)^{\textrm{\# 1s 
in top row}\,+\,\textrm{\# $\underline{0}$s in bottom row}}$. Of 
course this does not specify any element's color but no matter: that 
color will never be changed. For example, with $m=10$,
\[ 
\bordermatrix{
       &  1  &  2    &  3 & 4 & 5 & 6 &  7  & 8 & 9 & 10 \cr
  &       0  &  \uo  &  1 & 0 & 1 & 1 &  0  & 0 & 0 & 0 \cr
  &       1  &  1    &  1 & 0 & 1 & 1 & \uo & 0 & 0 & 0\cr
 }
\]
represents $K=\{3,5,6\},\ J=\{1,2,3,5,6\},\ A=\{2\},\ B=\{7\}$. Each 
column is one of $\ \genfrac{}{}{0pt}{}{0}{0}\ ,\ \genfrac{}{}{0pt}{}{0}{1},\ 
\genfrac{}{}{0pt}{}{1}{1},\ \genfrac{}{}{0pt}{}{\uo}{1},\ 
\genfrac{}{}{0pt}{}{0}{\uo}\ $ and the only further restriction is that 
there must be at least as many plain 0s terminating the bottom row as 
1s in the top row.

The weights of all configurations containing a colored vertex (= 
underlined matrix entry)  cancel out: Look for the leftmost underline 
and if its column is $\left(\begin{smallmatrix} 0 \\ \uo 
\end{smallmatrix}\right)$ then 
change that column to $\left(\begin{smallmatrix} \uo \\ 1 \end{smallmatrix}\right)$ 
and vice versa.
Now for a weight-reversing involution on 
matrix configurations with no underlines. The survivors are those 
with $K=\emptyset$ (i.e. top row all 0s) and $J$ a (possibly empty) 
initial segment of $[m]$ (i.e. bottom row is 1s followed by 0s). 
There are $m+1$ such, all of weight 1. The involution on the others 
is as follows. Look for the first $\left(\begin{smallmatrix} 1 \\ 1 
\end{smallmatrix}\right)$ column and consider $v$, the list of 
entries in the bottom row following this $\left(\begin{smallmatrix} 1 \\ 1 
\end{smallmatrix}\right)$; $v=$ all of bottom row if there is no $\left(\begin{smallmatrix} 1 \\ 1 
\end{smallmatrix}\right)$ column and note that, in this case, $v$ 
must contain a 01 string since $J$ is not an initial segment of $[m]$. If $v$ 
has no 01 string (and hence $v$ has the form $1^{i}0^{j}$ with $i\ge 
0$), replace the $\left(\begin{smallmatrix} 1 \\ 1 
\end{smallmatrix}\right)$ by $\left(\begin{smallmatrix} 0 \\ 0 
\end{smallmatrix}\right)$ and change the first 0 in $v$ to 1. There 
must be such a 0 since $\max\{J\} \le m- \vert K \vert$. On the other 
hand, if $v$ has a 01 string, change the column of the 0 in $v$'s 
last 01 string to $\left(\begin{smallmatrix} 1 \\ 1 
\end{smallmatrix}\right)$ and change the last 1 in $v$ to 0. An 
example is illustrated.

\begin{eqnarray} \label{bij1}
 \bordermatrix{
       &  1  &  2    &  3 & 4 & 5 & 6 &  7  & 8 & 9  \cr
  &       0  &  1    &  0 & 0 & 0 & 0 &  0  & 0 & 0  \cr
  &       1  &  1    &  0 & 1 & 0 & 1 &  1  & 1 & 0 \cr
  } & \longleftrightarrow &
 \bordermatrix{
       &  1  &  2    &  3 & 4 & 5 & 6 &  7  & 8 & 9  \cr
  &       0  &  1    &  0 & 0 & 1 & 0 &  0  & 0 & 0  \cr
  &       1  &  1    &  0 & 1 & 1 & 1 &  1  & 0 & 0 \cr
 }   \\ \nonumber
\substack{ \textrm{$v=0\,1\, 0\, 1\, 1\, 1\, 0$ contains a 01 
string\hspace*{5mm}} \\  
\textrm{The 0 of the last 01 in $v$ is in column 5\hspace*{5mm}} \\
\textrm{\vphantom{$\left(\begin{smallmatrix} 1 \\ 1 
\end{smallmatrix}\right)$}$v$'s last 1 is in column 8\hspace*{5mm}} }
  & & 
\substack{\hspace*{10mm} \textrm{$v= 1\, 1\, 0 \, 0$ contains no 01 string}\\
\textrm{\hspace*{10mm}The last $\left(\begin{smallmatrix} 1 \\ 1 
\end{smallmatrix}\right)$  is in column 5} \\
\textrm{\hspace*{10mm}$v$'s first 0 is in column 8} }
   \end{eqnarray}

Finally, we do the case $0\le k < m$:
\begin{multline*}
    \sum_{0 \le j \le i \le  m}^{}(-1)^{m-i}\binom{i}{j}
    \binom{m-k+j}{m-i}(1+q)^{i+j-k}(1-q)^{i-j}\\
    =\ \big( (m-k)+(m+1)q\big)q^{m-k-1}\hspace*{20mm}
\end{multline*}

Here we need to consider the $(m-k)$-set $E=[m+1,2m-k]$ as well as $M=[m]$. The 
5-tuple configurations are $(i,J,K,A,B)$ with $i \in 
[0,m],\ J$ a $j$-subset of $I:=[i],\ K$ an $(m-i)$-subset of $J\cup E,\ A$ 
an arbitrary $q$-colored subset of $(J \cup E)\backslash K\ $ and $B$ 
an arbitrary $q$-colored subset of $I\backslash J$, and weight 
$=(-1)^{\vert K \vert+\vert B \vert}$. Again we can drop the ``$i$'' 
by imposing the condition $\max\{J,B\}\le m -\vert K \vert$ and we 
will work with the equivalent matrix formulation which is the same as 
before except now the first row has an extra $m-k$ entries. 
Specifically, we have a $2 \times m$ main matrix each column being 
$\ \genfrac{}{}{0pt}{}{0}{0}\ ,\ \genfrac{}{}{0pt}{}{0}{1},\ 
\genfrac{}{}{0pt}{}{1}{1},\ \genfrac{}{}{0pt}{}{\uo}{1},\ $ or
$\ \genfrac{}{}{0pt}{}{0}{\uo}\ $ as before, an $(m-k)$-vector 
(considered as an extension of the top row of the matrix) each entry 
being $0,1$ or $\uo$, and the restriction that there are at least as 
many plain 0s terminating the bottom row as 
there are 1s in the top row and vector combined.
We have to expect complications due to the relatively complicated 
right hand side, $\big( (m-k)+(m+1)q\big)q^{m-k-1}$, and 
what we'll do is kill off weights in four steps, pruning our class of 
configurations at each step except for some easily counted survivors, each of weight 
1.

First, our earlier $\left(\begin{smallmatrix} 0 \\ \uo 
\end{smallmatrix}\right) \leftrightarrow \left(\begin{smallmatrix} \uo \\ 1 \end{smallmatrix}\right)$ 
involution kills off underlines in the main matrix. Henceforth only 
the $(m-k)$-vector, call it $v$, can contain underlines. 
Second, the bijection of 
(\ref{bij1}) kills all configurations where either the main matrix contains 
either $\left(\begin{smallmatrix} 1 \\ 1 
\end{smallmatrix}\right)$ or its bottom row contains 01 (or both).

We have now pruned our configurations to those with a 0-1 (no 
underlines) matrix with top row all 0s, and all 1s (if any) in the 
bottom row flush left. Those for which $v$ is all 0s are 
survivors---$(m+1)q^{m-k}$ of them. Otherwise, $v$ has the form 
$\uo^{i}aw$ with $i\ge 0,\ a=0$ or 1, and $w$ a (possibly empty) 
vector of 1s, 0s, $\uo$s.
Consider the configurations with bottom row all 0s and $a=0$. If 
$w$ in $\uo^{i}aw$ is all $\uo$s, we have a survivor---$(m-k)q^{m-k-1}$ of them. 
Otherwise, flipping the first 0/1 entry in $w$ ($0 \leftrightarrow 1$) is our 
involution. All that's left are configurations with (i) bottom row of 
the form $1^{i}0^{m-i}$ with $i\ge 1$ or (ii) $a=1$ (or both). Here 
the involution is: if $a=1$, change it to 0  and change the bottom row to 
$1^{i+1}0^{m-(i+1)}$ (there will be room for this extra 1 in the bottom row 
since $a=1$ forces a trailing 0 in the bottom row), and if $a=0$, change it to 1 
and change the bottom row to $1^{i-1}0^{m-(i-1)}$. We are done.

\end{document}